\documentclass[11pt]{amsart}

\usepackage[colorlinks=true, pdfstartview=FitV, linkcolor=blue, citecolor=blue, urlcolor=blue]{hyperref}

\usepackage{amsmath}
\usepackage{geometry}
\usepackage{layout}
\usepackage{graphicx}
\usepackage{amssymb}
\usepackage{epstopdf}
\usepackage{lscape}
\usepackage[utf8]{inputenc}
\usepackage{tikz}
\usepackage[shortlabels]{enumitem}

\definecolor{darkblue}{rgb}{0.0,0,0.7} % darkblue color
\definecolor{darkred}{rgb}{0.7,0,0} % darkred color
\definecolor{darkgreen}{rgb}{0, .6, 0} % darkgreen color

\newcommand{\defn}[1]{{\color{darkred}\emph{#1}}} % emphasis of a definition

\theoremstyle{definition}

\numberwithin{equation}{section}

\title{A module for the Delta conjecture}

\author{Mike Zabrocki}
\address[M. Zabrocki]{Department of Mathematics and Statistics,  York University, 4700 Keele Street, Toronto, 
Ontario M3J 1P3, Canada}
\email{zabrocki@mathstat.yorku.ca}
\urladdr{\url{http://garsia.math.yorku.ca/~zabrocki/}}
\begin{document}
%\layout

\begin{abstract}
We define a module that is an extension of the diagonal harmonics and
whose graded Frobenius characteristic is conjectured to be the symmetric function
expression which appears in `the Delta conjecture'
of Haglund, Remmel and Wilson \cite{HRW}.
\end{abstract}
\maketitle

%%%%%%%%%%%%%%%%%%%%%%%%%%%%%%%%%%%%%%%%%%%%%%%%%%
%\section{Introduction}
%%%%%%%%%%%%%%%%%%%%%%%%%%%%%%%%%%%%%%%%%%%%%%%%%%
The following account of the history of the conjecture stated here is necessarily abbreviated.
A number of recent expository articles \cite{vW, Hicks} cover the relevant
history of the shuffle conjecture and diagonal harmonics in greater detail. 

In \cite{Hai1}, Haiman defined the module of \defn{diagonal coinvariants} and conjectured
that it had dimension $(n+1)^{n-1}$.  The statement of the graded Frobenius image in terms of
Macdonald symmetric functions of this module first appears in a paper by
Garsia and Haiman \cite{GarHai} and an expression in terms of operators $\nabla$ and
$\Delta_f$ first appears in a paper by Bergeron, Garsia, Haiman and Tesler \cite{BGHT}.
In 2002, Haiman \cite{Hai2} published a proof that
$\nabla(e_n)$ was the formula for the $qt$-graded Frobenius image of the diagonal coinvariants.
In 2005, Haglund, Haiman, Loehr, Remmel and Ulyanov \cite{HHLRU} published
a conjectured combinatorial formula for the monomial expansion of $\nabla(e_n)$ 
and the equality became known as \defn{the shuffle conjecture}.
Carlsson and Mellit \cite{CM} recently published a proof of this conjecture.

Shortly after the proof of the shuffle conjecture was first posted,
Haglund, Remmel and Wilson \cite{HRW} defined \defn{the Delta conjecture} as the
equality of a combinatorial expression and the symmetric function  
$\Delta_{e_k}'(e_n)$. 
The expression for $k=n-1$ is special case of this formula that reduces to
the shuffle conjecture.  In January 2019, a workshop was held at BIRS in Banff titled
\textit{Representation Theory Connections to $(q,t)$-Combinatorics} \cite{Banff}.
At that meeting the author proposed the following conjecture as a
representation theoretic model for
the symmetric function expression $\sum_{k=1}^{n} z^{k-1} \Delta_{e_{n-k}}'(e_n)$.
A proof of this conjecture would
imply that the symmetric function expression is Schur positive.

%\section{The conjecture}
Fix an integer $n$ and let
$R_n:={\mathbb Q}[x_1, x_2, \ldots, x_n, y_1, y_2, \ldots, y_n, \theta_1, \theta_2, \ldots, \theta_n]$
be a polynomial ring in three $3n$ variables.  The $x_i$ and $y_i$ variables all commute
and commute with the $\theta_i$ variables, however the $\theta_i$ variables are Grassmannian,
that is, $\theta_i^2 = 0$ and $\theta_i \theta_j = - \theta_j \theta_i$ for $1 \leq i \neq j \leq n$.  Let
$I_n$ be the ideal of $R_n$ generated by elements
$$p_{r,s} := x_1^r y_1^s + x_2^r y_2^s + \cdots + x_n^ry_n^s\hbox{ for integers }0 < r+s\leq n
\qquad\hbox{ and }$$
$${\tilde p}_{r',s'} := x_1^{r'} y_1^{s'} \theta_1 + x_2^{r'} y_2^{s'} \theta_2
+ \cdots + x_n^{r'}y_n^{s'} \theta_n\hbox{ for integers }0 \leq r'+s' <n~.$$
The quotient $M_n := R_n/I_n$ will be referred to as the \defn{super-diagonal coinvariants}
(borrowing the `super' in the name from symmetric functions in superspace \cite{DLM}).
A consequence of Theorem 4.5 of \cite{OZ} implies the $p_{r,s}$ and ${\tilde p}_{r',s'}$ which
are generators of the ideal $I_n$ also algebraically generate the ring of $S_n$-invariants
of $R_n$.

The symmetric group acts on this ring by simultaneously permuting the indices of three sets of
variables.  The ring $R_n$ and the quotient $M_n$
are tri-graded by the degree in the three sets of variables.
For non-negative integers $a$, $b$ and $c$, let $M_n^{(a,b,c)}$ denote the homogeneous subspace of
$M_n$ of degree $a$, $b$ and $c$ in the respective variables $x_i$, $y_i$ and $\theta_i$
(for $1 \leq i \leq n$).  These
subspaces are symmetric group modules under this action of permutation of
the indices of the variables.

For a partition $\mu$ of $n$, let $\chi_{M_n^{(a,b,c)}}(\mu)$
denote the value of the character of the action of a permutation of cycle type
$\mu$ on the module $M_n^{(a,b,c)}$.
Define a \defn{$qtz$-graded Frobenius image} for super-diagonal coinvariants as
$${\mathcal F}_{qtz}( M_n ) = \sum_{a,b,c \geq 0}
\sum_{\mu} q^a t^b z^c \chi_{M_n^{(a,b,c)}}(\mu) \frac{p_{\mu}}{z_{\mu}}$$
where $p_\mu$ is the power sum symmetric function and $z_\mu$ is the expression
$\prod_{i=1}^{\mu_1} m_i! i^{m_i}$ with $m_i$ equal to the number of
parts of size $i$ in $\mu$.

The reader is referred to \cite{HRW} for the relevant background and notation
for the definition of
$\Delta'_f$.  It is defined as the symmetric function operator such that
$\Delta'_f( {\tilde H}_\mu[X;q,t]) = f[B_\mu(q,t)-1] {\tilde H}_\mu[X;q,t]$.
The purpose of this document is to state the following:\\

\noindent
{\bf Conjecture:}~~
For $n \geq 1$,
$${\mathcal F}_{qtz}( M_n ) = \Delta'_{e_{n-1} + z e_{n-2} + \cdots + z^{n-1}}(e_n).$$\vskip .1in

The special case of $z=0$ reduces to the theorem of Haiman for the diagonal coinvariants.
The idea for adding the anti-commuting set of variables to
this conjecture arose from a project first
proposed through the algebraic combinatorics seminar
at the Fields Institute \cite{Fields} in 2018--2019.  The project corresponds
to the $t=0$ case of this conjecture and we
referred to it as the \defn{super-coinvariants}.
When we later realized that the super-coinvariant modules are related to the $t=0$ case of the 
Delta conjecture, the extension to the super-diagonal coinvariants was a natural candidate for the
full Delta conjecture.

The Hall-Littlewood expansion of the $t=0$ case of the Delta conjecture from
\cite{GHRY,HRS} was used to arrive at this conjecture.
It is interesting to note other module structures have been proven
by Haglund, Rhoades, Shimozono and Wilson \cite{HRS, RW} to have
Frobenius image equal (up to some transformation) to $\Delta'_{e_{k}}(e_n)|_{t=0}$.

A program that uses Macaulay 2 and
Sage for computing the Frobenius image ${\mathcal F}_{qtz}( M_n )$ is available \cite{Zab}.
This program was used to verify the conjecture up to $n=6$.

\thanks{The author would like to thank the participants in the Fields seminar
\cite{Fields} including Nantel Bergeron, Laura Colmenarejo, John Machacek, Robin Sulzgruber
and Shu Xio Li.
The author would also like to thank Adriano Garsia for numerous conversations about these modules.}

\end{document}